\documentclass{amsart}
\usepackage{graphicx, color} 
\usepackage[margin=0.5cm,labelformat=parens]{subcaption} 
\usepackage{overpic}

\newcommand{\HH}{{\mathbb{H}}}
\newcommand{\RR}{{\mathbb{R}}}
\newcommand{\ZZ}{{\mathbb{Z}}}

\newcommand{\CC}{{\mathbb{C}}}

\newcommand{\calP}{{\mathcal{P}}}
\newcommand{\calD}{{\mathcal{D}}}
\newcommand{\wh}[1]{\widehat{#1}}
\newcommand{\bound}{{\partial}}
\DeclareMathOperator{\Vol}{Vol}
\DeclareMathOperator{\Li}{Li}

\newtheorem{theorem}{Theorem}
\title{Volumes of ideal hyperbolic drums}
\author[E.~Chesebro]{Eric Chesebro}
\author[J.~Lanier]{Justin Lanier}
\author[E.~McQuire]{Emma N. McQuire}
\author[J.~Morgan]{James Morgan}
\author[J.~Purcell]{Jessica S. Purcell}
\author[H.~Segerman]{Henry Segerman}

\address{Eric Chesebro \newline Department of Mathematics, University of Montana\newline Missoula, MT, 59812, USA \\  eric.chesebro@umontana.edu}
\address{Justin Lanier \newline Department of Mathematics, Louisiana State University\newline Baton Rouge, LA, 70802, USA \\  justin.lanier@lsu.edu}
\address{Emma N. McQuire \newline School of Mathematics, Monash University \newline Clayton, VIC 3800, Australia \\
emma.mcquire@monash.edu}
\address{James Morgan \newline School of Mathematics and Statistics, University of Sydney, \newline Sydney, NSW 2006, Australia \\
james.morgan@sydney.edu.au}
\address{Jessica S. Purcell \newline School of Mathematics, Monash University \newline Clayton, VIC 3800, Australia \\
jessica.purcell@monash.edu}
\address{Henry Segerman \newline Department of Mathematics, Oklahoma State University\newline Stillwater, OK, 74078, USA \\ henry.segerman@okstate.edu}

\begin{document}
\begin{abstract}
Milnor computed the volumes of ideal hyperbolic prisms as part of an effort to construct 3-manifolds whose volumes are finite rational sums of the Lobachevsky function evaluated at rational multiples of $\pi$. Motivated by these results and with an eye to related applications, we prove a volume formula for arbitrary ideal hyperbolic antiprisms, also called drums.
\end{abstract}
\maketitle
\begin{figure}[h]
\includegraphics[width=0.5\textwidth]{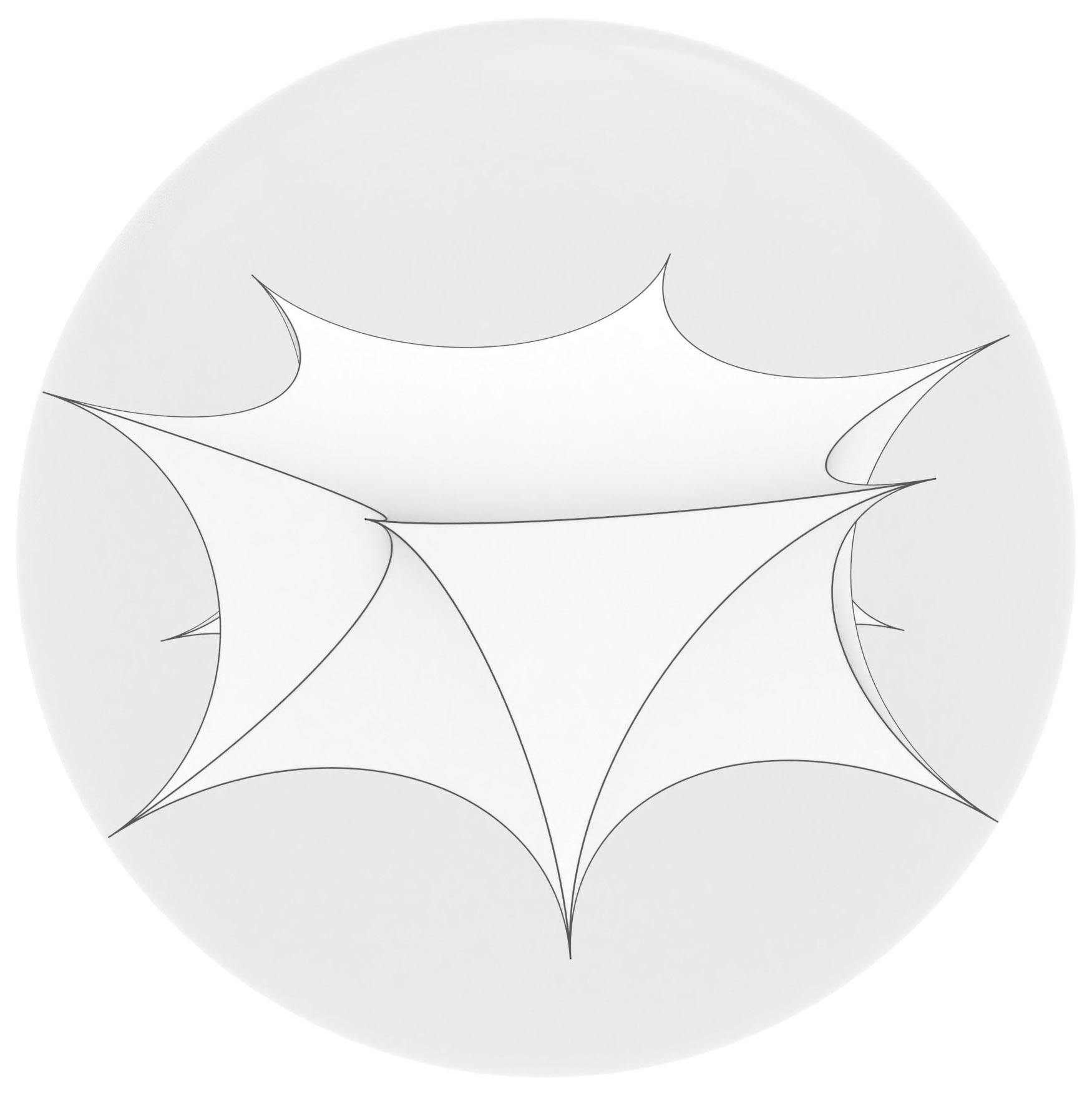}
\caption{A $6$-drum drawn in the Poincar\'e ball model of $\HH^3$. This drum corresponds to the parameters $n=6$, translation length $0.78$, and rotation angle $2\pi/9$.}
\label{Fig:Drum}
\end{figure}

\section{Introduction}

Fix $n \in \ZZ_{\geq 3}$ and let $\calP$ be a totally geodesic regular ideal $n$-gon in hyperbolic 3-space $\HH^3$.  Then there is an elliptic isometry $\sigma$ which rotates through an angle $2\pi/n$ about an axis $\alpha$ and restricts to a rotational symmetry on $\calP$.  Let $\varphi$ be a loxodromic isometry with axis $\alpha$, a non-trivial translation length, and a (possibly trivial) rotation angle.  We refer to the convex hull of $\calP \cup \varphi(\calP)$ as an {\em $n$-drum}. Combinatorially, an $n$-drum is an antiprism, excepting when it is a prism when the rotation angle is a multiple of $2\pi/n$. The purpose of this paper is to give a general formula for the volume of an $n$-drum.

In Chapter 7 of \cite{Th_notes}, Milnor produces a formula for the volumes of $n$-drums with rotation angle 0. He does so in pursuit of hyperbolic polyhedra where the group generated by the reflections in their sides form a discrete group of isometries of $\HH^3$. For these prism drums, he shows this condition is met for each $n$ for exactly two pairs of dihedral angle parameters: $(\pi/3,\pi/3)$ and $(\pi/4,\pi/2)$, the first entry giving the dihedral angle along the rims of the drum and the second the dihedral angle along the drum's vertical edges. In Chapter 6 of \cite{Th_notes}, Thurston states a formula for the volume of an antiprism $n$-drum with all dihedral angles $\pi/2$, mentioning that it can be derived using methods similar to those used by Milnor.

We are motivated to generalize Milnor's results with some particular applications in mind. Some important classes of hyperbolic 3-orbifolds have compact cores which can be decomposed into ideal polyhedral cells which include drums.  Several families of examples are given in Chapter 6 of \cite{Th_notes}; these include chain link complements and certain symmetric (closed) braid complements.  More recently, it is shown in \cite{chesebro_geom} and \cite{CEP} that the convex cores of Heckoid orbifolds decompose into combinations of drums and Sakuma--Weeks triangulations.  Moreover, cyclic branched covers of 2-bridge link complements, branched over crossing circles, decompose into combinations of drums and Sakuma--Weeks triangulations \cite{CEP}.

The work in this paper was carried out at the workshop {\it Triangulations in Low-Dimensional Topology} at the  MATRIX Institute in Victoria, Australia. We thank MATRIX and the workshop organizers for their support, organizational work, and generous hospitality.

\section{Methods and result}

\begin{figure}[htbp]
\subfloat[Suspension of the drum.]{
\includegraphics[width=0.42\textwidth]{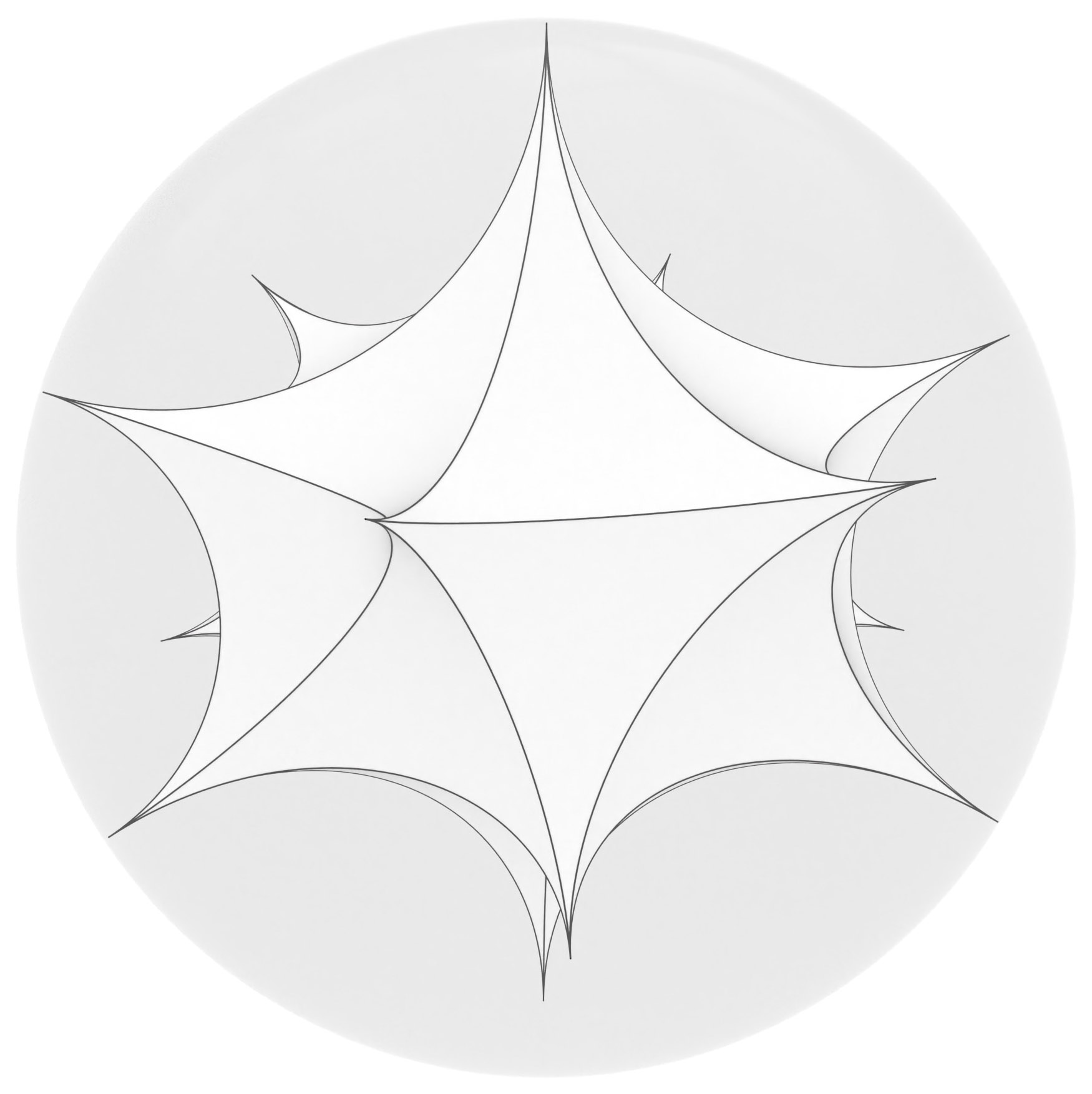}
\label{Fig:SuspensionDrum}
}
\subfloat[Remove two pyramids to recover the drum.]{
\includegraphics[width=0.42\textwidth]{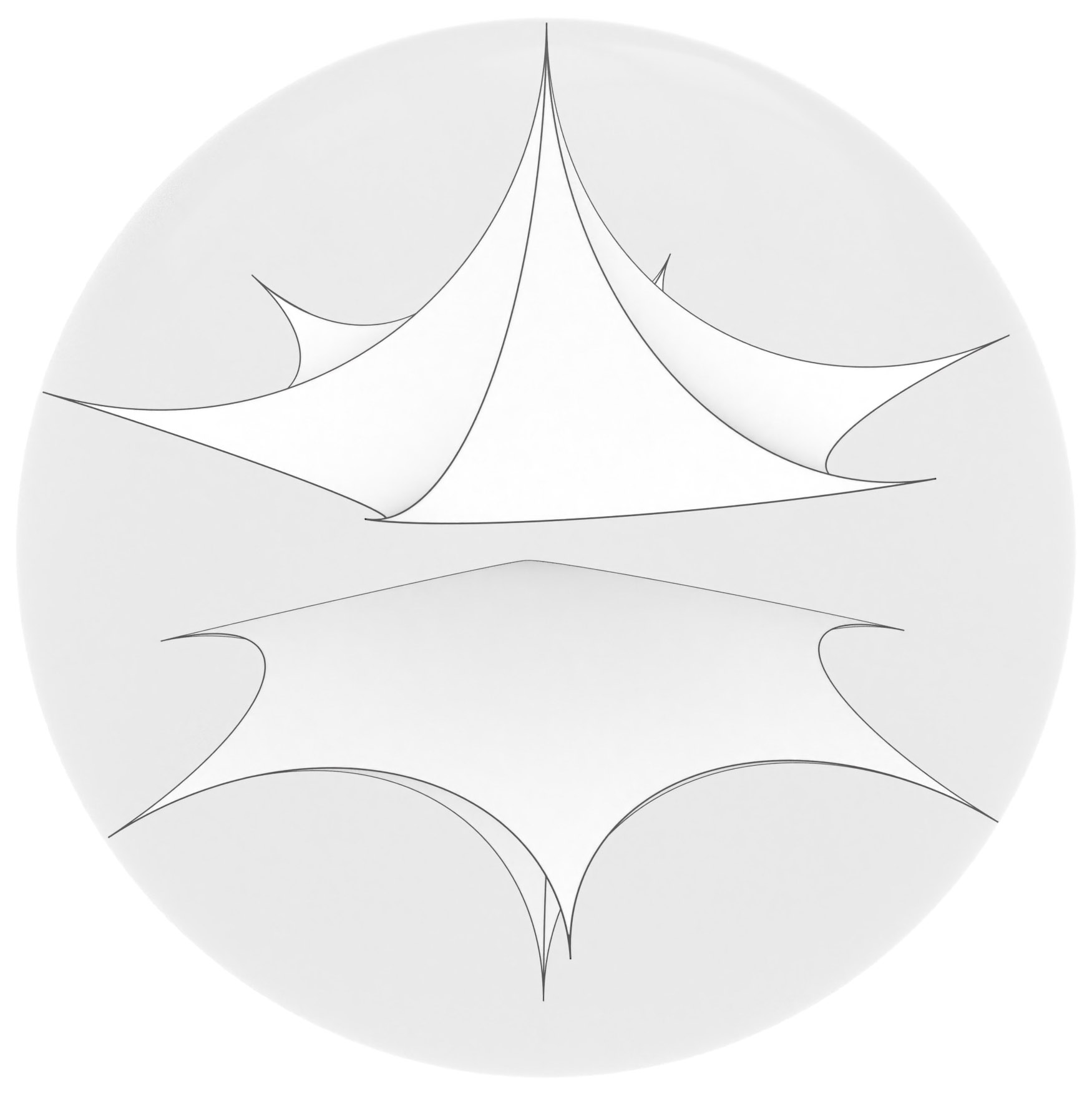}
\label{Fig:TwoPyramids}
}
\caption{The volume of the drum is the volume of the suspended drum minus twice the volume of an ideal pyramid.}
\label{Fig:Suspension}
\end{figure}

Suppose that $\calD \subset \HH^3$ is an $n$-drum, as described above.  We use the upper-half space model $\{z + tj \mid z\in \CC, t\in \RR^+\}$ of $\HH^3$ so that $\partial \HH^3$ is the Riemann sphere $\widehat{\CC} = \CC \cup \{ \infty \}$.   By applying an isometry to $\calD$, we may assume that the geodesic axis $\alpha$ travels from $0$ to $\infty$ and the ideal points of $\calP$ are 
\[\left\{ e^{2 k \pi i/n} \:\, \Big| \:\, k=0,1,2, \ldots, n-1\right\}\]
and the loxodromic isometry $\varphi$ acts on $\widehat{\CC}$ through multiplication by $re^{i \theta}$, where $|r|>1$ and $0\leq \theta < 2\pi/n$.  Refer to the polygons $\calP$ and $\varphi \calP$ as the {\em bottom and top heads} of the drum $\calD$.  Notice that the axis $\alpha$ meets the boundary of the drum at the centers of the heads at the points $j$ and $rj$.  The complement of the heads in $\bound \calD$ is called the {\it shell} and the geodesics in the shell that connect ideal points are called {\it laces}.  The heads meet the shell in the upper and lower {\it rims}; the rims coincide with $\bound \calP$ and $\bound \varphi \calP$.   Notice that, with this normalization, the $n$-drum $\calD$ is determined by the single complex number, $re^{i\theta}$, or two real parameters, $r$ and $\theta$.

It is tempting to triangulate a fundamental domain for the action of the elliptic $\sigma$ on $\calD$.  However, these triangulations will include tetrahedra with some finite vertices.  Since the volume formulas for ideal tetrahedra are nicer, we seek a different approach.

Define the {\it suspension} $\wh{\calD}$ of $\calD$ to be the convex hull of the union of $\calD$ with the axis of $\varphi$ (the geodesic between $0$ and $\infty$).  See Figure~\ref{Fig:SuspensionDrum}.  The complement of $\calD$ in $\wh{\calD}$ consists of a pair of regular ideal pyramids.  The upper pyramid is the cone of $\varphi \calP$ to $\infty$ and the lower pyramid is the cone of $\calP$ to $0$.  See Figure~\ref{Fig:TwoPyramids}.  Because the two pyramids are isometric,
\[ \Vol(\calD) \ =\ \Vol\left( \wh{\calD}\right) - 2 V\]
where $V$ is the volume of one pyramid.  If $\varphi^{-1}$ is applied to the upper pyramid, its union with the lower pyramid is the suspension $\wh{\calP}$, the convex hull of the union of $\calP$ with the axis of $\varphi$.  Both $\wh{\calP}$ and $\wh{\calD}$ admit nice decompositions into ideal tetrahedra and 
\[ \Vol(\calD) \ =\ \Vol\left( \wh{\calD}\right) - \Vol(\wh{\calP}).\]

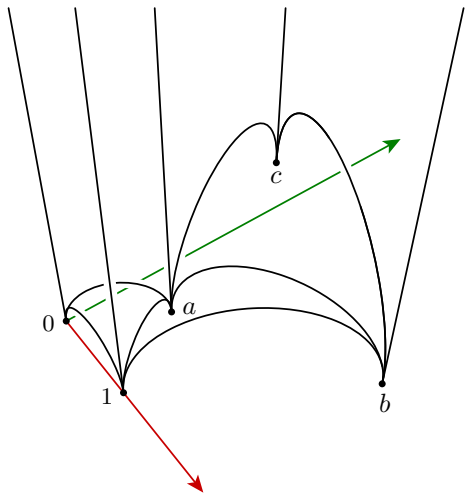
\begin{figure}[htbp] 
\begin{overpic}[width=0.48\textwidth, percent]{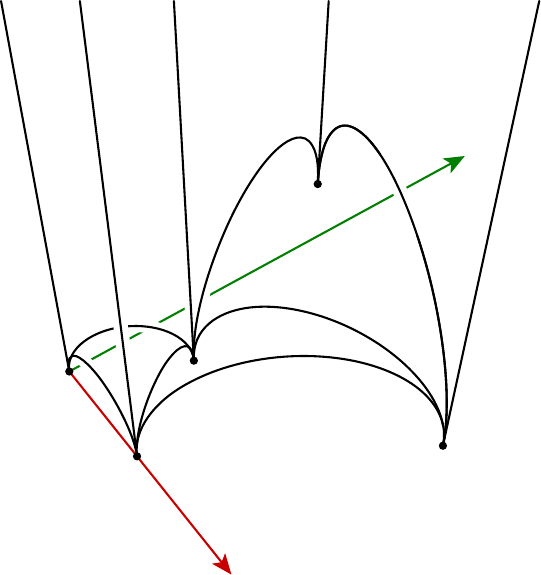}
\put(7.2,33.4){\small $0$}
\put(19.2,18.3){\small $1$}
\put(36,36.8){$a$} 
\put(54,63.7){$c$} 
\put(76.5,16.7){$b$}
\end{overpic}
\caption{Three ideal tetrahedra, drawn in the upper half space model.}
\label{Fig:UHS_ideal}
\end{figure}

Next, we describe a decomposition of $\wh{\calD}$ into ideal tetrahedra.  For notational convenience, take $\alpha=2\pi/n$ and write
\begin{align*}
a&= e^{i\alpha} & b&=re^{i\theta} & c&=re^{i(\alpha+\theta)}.
\end{align*}
Then $c=ab$.  Define three ideal tetrahedra by listing their ideal vertices
\begin{align*}
\tau_1&=[ 0,1,a,\infty] & \tau_2&=[ 1,a,b,\infty] & \tau_3&=[a,b,c,\infty].
\end{align*}
Shown in Figure~\ref{Fig:UHS_ideal}, the union of these tetrahedra form a fundamental domain for the action of $\langle \sigma \rangle$ on $\wh{\calD}$.  The tetrahedron $\tau_1$ is, by itself, a fundamental domain for the action of $\langle \sigma \rangle$ on $\wh{\calP}$.  Therefore,
\[ \frac{\Vol(\calD)}{n} \ =\  \Vol(\tau_2) + \Vol(\tau_3).\]

\begin{figure}[htbp]
\subfloat[The suspension $\wh{\calP}$ of $\calP$ is an ideal bipyramid.]{
\includegraphics[width=0.42\textwidth]{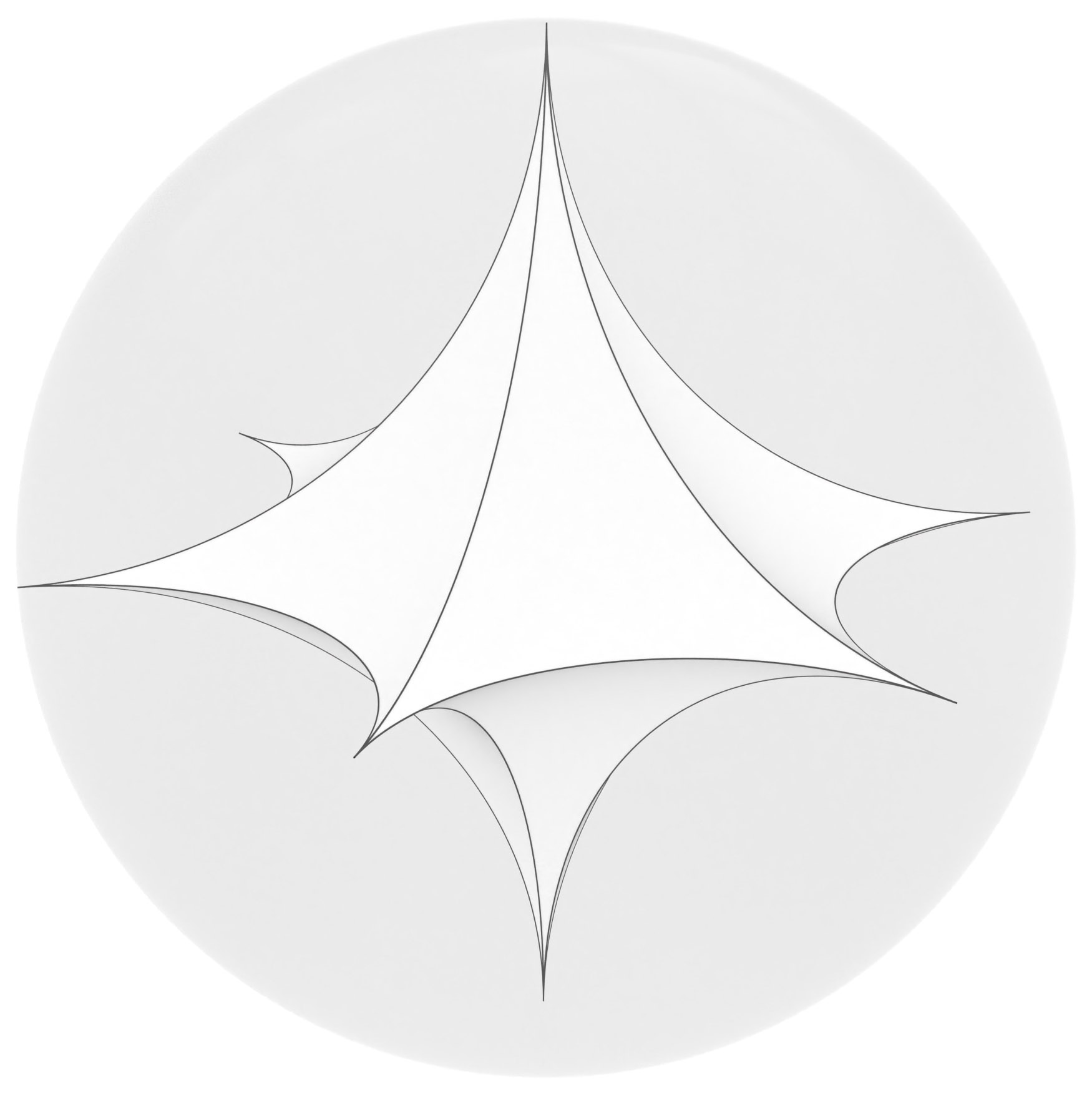}
\label{Fig:Bipyramid}
}
\subfloat[The tetrahedron $\tau_1$ is a fundamental domain for the rotational symmetry of $\wh{\calP}$.]{
\includegraphics[width=0.42\textwidth]{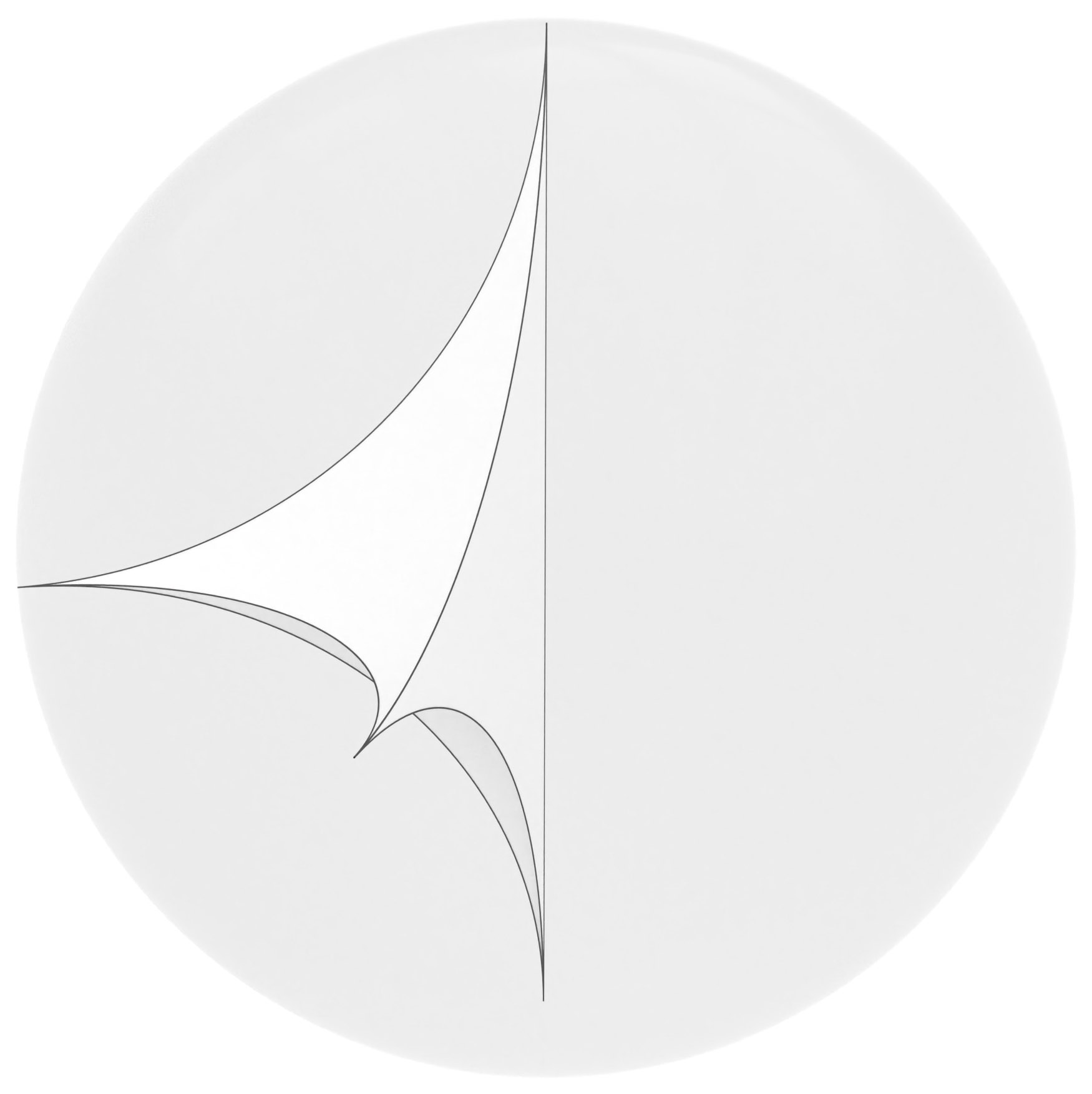}
\label{Fig:OneTetBipyramid}
}
\caption{The bipyramid $\wh{\calP}$ and its fundamental domain $\tau_1$.}
\end{figure}

The formulas of \cite{NZ} allow us to compute the volumes of $\tau_2$ and $\tau_3$ in terms of their respective edge parameters $z_2$ and $z_3$ for their common edge $E$,  running from $a$ to $\infty$.  A quick calculation shows
\begin{align}
\label{eq: 1} z_2&= \frac{b-a}{1-a} = \frac{e^{i\alpha}-re^{i\theta}}{e^{i\alpha}-1} 
\end{align}
and
\begin{align}
\label{eq: 2} z_3&= \frac{c-a}{b-a}= \frac{re^{i(\alpha+\theta)}-e^{i\alpha}}{re^{i\theta}-e^{i\alpha}}.
\end{align}
Using formula (38) in Section 5 of \cite{NZ}, the volume of $\tau_k$ is
\begin{align} \label{eq: 3}  \Vol(\tau_k) &= \Im \left(\Li_2(z_k)\right) +\ln|z_k|\arg(1-z_k)\end{align}
where 
\[ \Li_2(z) \ = \ \sum_{m=1}^\infty \frac{z^m}{m^2}\]
is the dilogarithm function.  This allows us to compute the volume of the drum $\calD$:

\begin{theorem}
The volume of a hyperbolic $n$-drum $\calD$ determined by parameters $r$ and $\theta$ is
\[ \Vol(\calD) \ =\  n\left(\Vol(\tau_2) + \Vol(\tau_3)\right),\]
where the volumes $\Vol(\tau_2)$ and $\Vol(\tau_3)$ may be computed using formulas (\ref{eq: 1}), (\ref{eq: 2}), and (\ref{eq: 3}).
\end{theorem}

\textit{Matching up with Milnor's result.} For $n \geq 3$, Milnor computes the volumes of ideal hyperbolic $n$-prisms in terms of their two dihedral angle parameters: $\alpha$ for the dihedral angles about the rim edges and $\beta$ for the dihedral angles about the vertical edges. These parameters have the linear relation $\beta=\pi-2\alpha$, and so we may parametrize these $n$-prisms by $n$ and $\alpha$ and denote them $\mathcal{N}_n(\alpha)$. These $n$-prisms are exactly $n$-drums with $\theta=0$ but parametrized by their rim dihedral angles instead of their heights. Here is Milnor's volume formula for $\mathcal{N}_n(\alpha)$ after making the substitution for $\beta$ \cite[Ch 7, Eq (6)]{Th_notes}:
\[
\Vol(\mathcal{N}_n(\alpha))=n(2L((\pi-2\alpha))/2)+L(\alpha+\pi/n)+L(\alpha-\pi/n))
\]

In the formula, $L$ is the Lobachevsky function:

\[
L(x)=- \int_{0}^{x}\log|2\sin t| \ dt.
\]

By applying an angle of parallelism formula and a formula for the summit angles of a Saccheri quadrilateral in terms of its base and legs, we have the following function for $\alpha$ in terms of $d$ and $n$:

\[
\alpha(d,n)=\pi/2-\arctan(\sinh(-\log(\tan(\pi/(2n)))\cdot\tanh(d/2))
\]

Plugging in this expression for $\alpha$ in Milnor's volume formula yields a formula for the $n$-prism's volume in terms of $n$ and $d$. Though they are quite different in form, this formula agrees with our general formula for the volume of an $n$-drum for any fixed $n$ with $\theta=0$ and $r=e^d$.
\section{Computations and open problems}

Consider the parameter spaces for $n$-drums $D_n(d, \theta)$ where $n$ is the order of rotational symmetry of the drum, $d$ is the translation length of the loxodromic $\varphi$, and $\theta$ is its rotation angle. The appropriate domain to have unique representatives of congruence classes of $n$-drums is the domain 
\[ \left\{(d,\theta) \in \mathbb{R}^2 \:\, \Big| \:\,  0<d,\ 0\leq \theta \leq \pi/n\right\}.
\]
Again, we have $d(\varphi)=\log(r(\varphi))$, where $r$ is the parameter from our upper-half space model of an $n$-drum. We conclude with some observations about the volume functional on this parameter space, which we invite others to investigate further.

For each fixed $n$, there appears to be a unique global maximum of $\Vol(D_n(d, \theta))$. We give a table of approximate maximizing parameters in Table~\ref{table:1}.

For each $n$, the global maximum appears to occur when $\theta=\pi/n$, the ``half-click" regular antiprism. This would be conveniently explained if volume were also maximized for fixed $n$ and for every $d$ at $\theta=\pi/n$, but this is not quite the case. For fixed $n$ and fixed $d$ greater than a smallish threshold (for instance, about $0.32$ for $n=3$), $\Vol(D_n(d, \theta))$ appears to be increasing and unimodal, with a single maximum at $\theta=\pi/n$. But for each $n$ and for sufficiently small values of $d>0$, $\theta=\pi/n$ is instead a local minimum for volume.
 
For large $n$, the increment between successive maximal volumes appears to approach twice the volume of a regular ideal tetrahedron, approximately $2.02988$. Another way to phrase the same observation is that, as $n$ tends to $\infty$, the ratio of the maximal volume of an $n$-drum to $n$ appears to approach twice the volume of a regular ideal tetrahedron.

Finally, for fixed $n$ and arbitrary $\theta$, it appears that as $d$ tends to $\infty$, 
$\Vol(D_n(d, \theta))$ tends to $\Vol(\wh{\calP})$. An informal justification for this is that as $d$ tends to $\infty$, $D_n(d, \theta)$ tends toward being two ideal $n$-pyramids connected at their apex cusps as the body of the drum becomes an increasingly long and narrow polygonal tube.
\begin{center}
\begin{table}
\caption{For various $n$, the approximate maximal volume of an $n$-drum and the approximate $d$ and $\theta$ values that yield this maximum.}
\begin{tabular}{||r l l l||} 
 \hline
 $n$ & Max vol & $d$ & $\theta$ \\ [0.5ex] 
 \hline\hline
 3 & 3.66386 & 1.31696 & 1.0472 \\
 \hline
 4&6.17545&1.12838&0.785398 \\
 \hline
 5&8.54285&0.962424&0.628319 \\
 \hline
 6&10.8158&0.831443&0.523599 \\
 \hline
  100&202.903&0.0543961&0.0314159 \\
 \hline
 10000&20298.8&0.00054414&0.000314157 \\[.5ex] 
 \hline
\end{tabular}
\label{table:1}
\end{table}
\end{center}
\bibliographystyle{amsplain}
\bibliography{biblioD}

\providecommand{\bysame}{\leavevmode\hbox to3em{\hrulefill}\thinspace}
\providecommand{\MR}{\relax\ifhmode\unskip\space\fi MR }
\providecommand{\MRhref}[2]{%
  \href{http://www.ams.org/mathscinet-getitem?mr=#1}{#2}
}
\providecommand{\href}[2]{#2}
\begin{thebibliography}{1}

\bibitem{chesebro_geom}
E.~Chesebro, \emph{Geometry for {K}leinian groups generated by a parabolic
  pair}, arXiv:2512.17044 (2026).

\bibitem{CEP}
E.~Chesebro, A.~Elzenaar, and J.~Purcell, \emph{Forthcoming work}.

\bibitem{NZ}
W.~D. Neumann and Don Zagier, \emph{Volumes of hyperbolic three-manifolds},
  Topology \textbf{24} (1985), no.~3, 307--332. \MR{815482}

\bibitem{Th_notes}
W.~P. Thurston, \emph{The geometry and topology of 3-manifolds}, mimeographed
  lecture notes, 1979.

\end{thebibliography}

\end{document}